\documentclass[12pt]{article}
\usepackage{amssymb}
\usepackage{amsmath}
\usepackage{algorithmic}
\usepackage{graphicx}
\usepackage{amsthm}
\theoremstyle{definition}
\theoremstyle{plain}
\newtheorem*{theorem}{\underline{Theorem}}
\newenvironment{algorithm}[1]{
  \begin{center}
    {\bf Algorithm: #1}\\*
    \begin{tabular}{|p{12.5cm}|} \hline
}
{
 \\ \hline
 \end{tabular}
 \end{center}
}
\newcommand{\Gr}{Gr\"obner }
\newcommand{\Z}{\mathbb{Z}}
\newcommand{\K}{\mathbb{K}}
\newcommand{\N}{\mathbb{N}}
\newcommand{\M}{\mathbb{M}}
\newcommand{\R}{\mathbb{R}}

\newcommand{\var}{\mathop{\mathrm{var}}\nolimits}
\newcommand{\ndg}{\mathop{\mathrm{ndg}}\nolimits}
\newcommand{\dg}{\mathop{\mathrm{dg}}\nolimits}
\newcommand{\nvr}{\mathop{\mathrm{nvr}}\nolimits}
\newcommand{\bnm}{\mathop{\mathrm{bnm}}\nolimits}
\newcommand{\lcm}{\mathop{\mathrm{lcm}}\nolimits}
\newcommand{\lm}{\mathop{\mathrm{lm}}\nolimits}

\newcommand{\anc}{\mathop{\mathrm{anc}}\nolimits}
\newcommand{\nmp}{\mathop{\mathrm{nmp}}\nolimits}

\newcommand{\bin}{\mathop{\mathrm{bin}}\nolimits}

\begin{document}

\title{\bf Janet Bases of Toric Ideals}
\author{Vladimir P. Gerdt \\
       Laboratory of Information Technologies\\
       Joint Institute for Nuclear Research\\
       141980 Dubna, Russia \\
       gerdt@jinr.ru
\and
       Yuri A. Blinkov \\
       Department of Mathematics and Mechanics \\
       Saratov State University \\
       410071 Saratov, Russia \\
       blinkovua@info.sgu.ru}
\date{}
\maketitle
\begin{abstract}
In this paper we present a version of the general
polynomial involutive algorithm for computing Janet bases
specialized to toric ideals. The relevant data structures are
Janet trees which provide a very fast search for a Janet divisor.
We broach also efficiency issues in view of application of the
algorithm presented to computation of toric ideals.
\end{abstract}

\section{Introduction}
We consider the problem of computing a Janet basis of a {\em toric}
ideal $\mathcal{I_A}$ in $\K[\mathbf{x}]\equiv \K[x_1,\ldots,x_n]$
generated by binomials of the form~\cite{Sturmfels96}
$$ \mathcal{I_A} = \{\ \mathbf{x^{u}-x^{v} \mid u,v}\in \N^n,\
 \pi(\mathbf{u})=\pi(\mathbf{v}),\ \gcd(\mathbf{x^{u}},\mathbf{x^{v}})=1\ \}\,. $$
Here $\mathbf{u,v}\in \N^n$ and $\pi$ is the semigroup homomorphism
$$ \pi\ :\ \N^n \rightarrow \Z^d,\quad \mathbf{u}=\{u_1,\ldots,u_n\}\rightarrow
u_1\mathbf{a}_1+\cdots +u_n\mathbf{a}_n$$
where $\mathbf{a}_i\in \Z^d$ $(1\leq i\leq n)$.

Given a set of binomials generating a toric ideals, the problem of constructing
its \Gr basis is usually (except small problems) rather expensive from the
computational point of view~\cite{BLR99}.  In practice, for this particular
problem, one typically deals with a large number $n$ of variables
and their degrees. If $d$ is the maximal degree of the initial binomials, then
the degree of a reduced \Gr basis is bounded by~\cite{KM96}
$$ 2\cdot \left(\frac{d^2}{2}+d\right)^{2^{n-1}}\,.$$
But for all that the reduced \Gr basis is also binomial since
the binomial structure is preserved during the Buchberger
algorithm~\cite{AL94,Pot94}. Similarly, the involutive algorithms~\cite{GB}
based on the sequential multiplicative reductions of nonmultiplicative
prolongations of the intermediate polynomials preserve the binomial structure.
The output involutive basis which is also a \Gr basis though
generally redundant.

Thus, unlike construction of reduced \Gr bases or involutive bases
for general polynomial ideals, the integer arithmetics which may
take most of computing time is not important for binomial ideals.
In this case a fast search of monomial divisors for performing
reductions of $S$-polynomials may become crucial in acceleration of
computations.

Recently~\cite{GBY1,GBY2} we designed and implemented involutive
algorithms specialized to constructing Janet bases of monomial and
polynomial ideals. Janet division as well as any other involutive
division~\cite{GB} provides uniqueness of an involutive divisor in
a polynomial set with co-prime leading monomials. This allows one
to organize a very fast search for a Janet divisor using special
data structures for intermediate polynomial sets called Janet
trees.

The main goal of this paper is to discuss the issue of practical
efficiency in computing Janet bases of toric ideals based on the use
of Janet trees. Since one of the most important applications of toric
ideals is integer programming we shortly describe this
application~\cite{CoTr91} in the next section.

\section{Toric Ideals and Integer Programming}
Let $ \mathcal{A}$ be a matrix of dimension $m \times n$ with
integer entries and $\mathbf{b} \in \Z^m$, $\mathbf{c} \in \Z^n$
be vectors. The following optimization problem
$$
 \min\{\ {\mathbf{c}^T \mathbf{x} \mid \mathbf{x} \in \N^n,\,
 \mathcal{A}\mathbf{x}=\mathbf{b}\ }\}
$$
is called a problem of {\em integer programming}.

We shall assume that $\mathbf{c}\in \N^n$. If there exists vector
$\mathbf{x}_0$ satisfying $ \mathcal{A}\mathbf{x}_0=\mathbf{b}$,
$\mathbf{x}_0\in \N^n$, then the problem of finding minimum of
function $\mathbf{c}^T \mathbf{x}$ can be reduced to all kinds of
transformation of the initial vector state $\mathbf{x}_0$ using
$\ker(\mathcal{A})$.

The problem of determining $\ker( \mathcal{A})$ can be formulated
in terms of toric ideals. Indeed, every vector $\mathbf{u} \in
\ker(\mathcal{A})$ may be uniquely represented as
$\mathbf{u}=\mathbf{u}^{+}-\mathbf{u}^{-}$ where both
$\mathbf{u}^{+}$ and $\mathbf{u}^{-}$ are nonnegative and have
disjoint support.  Associate symbol $v_i$ with the $i$-th column
of matrix $\mathcal{A}$. Then the ideal
$$
 \mathcal{I_A}=\{\ \mathbf{v}^{u^+}-\mathbf{v}^{u^-}\mid \mathbf{u}^{+}-\mathbf{u}^{-}=
 \mathbf{u} \in \ker( \mathcal{A})\ \}
$$
associated with $\ker(\mathcal{A})$ is toric. Given the initial solution
$\mathbf{x}_0$, the optimal solution can be found as
follows~\cite{CoTr91}:
\begin{enumerate}
  \item Construct a basis of the toric ideal $ \mathcal{I_A}$.
  \item Construct a reduced \Gr basis or an involutive basis of $ \mathcal{I_A}$ with
    respect to the admissible monomial ordering $\succ_\mathbf{c}$ generated by
    vector $\mathbf{c}$.
  \item Reduce monomial $\mathbf{v}^{\mathbf{x}_0}$ modulo the constructed basis to obtain
  the optimal solution.
\end{enumerate}
Therefore, the reduced \Gr basis or any involutive basis of the associated toric ideal
$ \mathcal{I_A}$ provide an algorithmic tool for solving the problem of integer
programming.

\section{Janet Bases of Toric Ideals}
\subsection{Definition of Janet Basis}
In our papers~\cite{GB} the \Gr bases of special type, called {\em
involutive} and based on the concept of {\em involutive division}
were introduced. Given a set of coprime monomials and an
involutive division, any monomial may have at most one involutive
divisor in the set. This property of the involutive division
allows one to design an efficient search for the involutive
divisor using the method of separative monomials~\cite{B01} for a
general involutive division or Janet trees~\cite{GBY1} for Janet
division.

Because of a larger number of variables and unimportance of
integer arithmetical operations over coefficients of the
binomials, the practical complexity of an algorithm for
construction of \Gr or Janet bases is caused by an enormous number
of binomials arising in computation of the basis. A faster search for
divisors may accelerates the computation substantially.

By definition of Janet division~\cite{GB} (which formalizes the pioneering ideas
of Janet~\cite{Janet}) induced by the order
\begin{equation}
x_1 \succ x_2 \succ \ldots \succ x_n \label{order}
\end{equation}
on $\mathbf{x}$, a
polynomial set $F$ is partitioned into the groups
labeled by non-negative integers
$d_1,\ldots,d_i$:
$$
[d_1,\ldots,d_i]=\{\ f\ \in F\ |\ d_j=\deg_j(\lm(f)),\ 1\leq j\leq i\ \}
$$
where $\deg_i(u)$ denotes the degree of $x_i$ in monomial $u$ and $\lm(f)$
denotes the leading monomial of $f$.
A variable $x_i$ is called {\em (Janet) multiplicative}
for $f\in F$ if $i=1$ and
$$\deg_1(\lm(f))=\max\{\ \deg_1(\lm(g))\ |\ g\in F\ \},$$
or if
$i>1$, $f\in [d_1,\ldots,d_{i-1}]$ and
$$\deg_i(\lm(f))=\max\{\ \deg_i(\lm(g))\ |\ g\in
[d_1,\ldots,d_{i-1}]\ \}\,.$$
If a variable is not multiplicative
for $f\in F$, it is called {\em (Janet) nonmultiplicative} for $f$. In the
latter case we shall write $x_i\in NM_J(f,F)$. $u\in \lm(F)$ is a {\em Janet divisor}
of $w\in \M$, if $u\mid w$ and monomial $w/u$ contains only multiplicative
variables for $u$. In this case we write $u\mid_J w$.

Let $\lm(F)=\{\ \lm(f) \mid f\in F\ \}$. Then a polynomial set $F$ is called
{\em Janet autoreduced} if each term in every $f\in F$ has no Janet divisors
among $\lm(F)\setminus \lm(f)$. A polynomial $h$ is said to be in the
{\em Janet normal form modulo $F$} if every term in $h$ has no Janet divisors
in $\lm(F)$. In that follows $NF_J(f,F)$ denotes the Janet normal form $f$
modulo $F$.

A Janet autoreduced set $F$ is called a {\em Janet basis} if
\begin{equation}
 (\forall f\in F)\ (\forall x\in NM_J(f,F))\ \ [\ NF_J(f\cdot x,F)=0\ ]\,. \label{J_basis}
\end{equation}
A Janet basis $G$ is called {\em minimal} if for
any other Janet basis $F$ of the same ideal the inclusion
$lm(G)\subseteq \lm(F)$ holds. If both $G$ and $F$ are monic this inclusion
implies $G\subseteq F$. A Janet basis is a \Gr one, though generally not
reduced. However, similarly to a reduced \Gr basis, a monic minimal Janet
basis is uniquely defined by an ideal and a monomial order.
In that follows we deal with minimal Janet bases only and omit the word "minimal".

\subsection{Janet Trees and Search for Janet Divisor}
Consider now a binary Janet tree~\cite{GBY1}
whose structure reflects the above partition of elements in $U$ into the groups
which sorted in the degrees of variables within every group. Before description of
the general structure of Janet trees we explain it in terms of the
concrete example~\cite{GBY1}
$$U=\{x^2y,xz,y^2,yz,z^2\},\qquad (x\succ y\succ z)$$
and portray it in the form of Janet tree as shown below. In doing so, the monomials
in set $U$ are assigned to the leaves of the tree. The monomial with
increased by one degree of the current variable is assigned to the left child
whereas the right child points at the next variable with respect to chosen ordering.
In contrast to Janet tree presented in paper~\cite{GBY1}, the below tree takes into
account sparseness of monomials that is inherent in toric ideals. The related
information is given in pairs of integers placed in brackets where the first element
represents the number of current variable and the second one represents its degree.

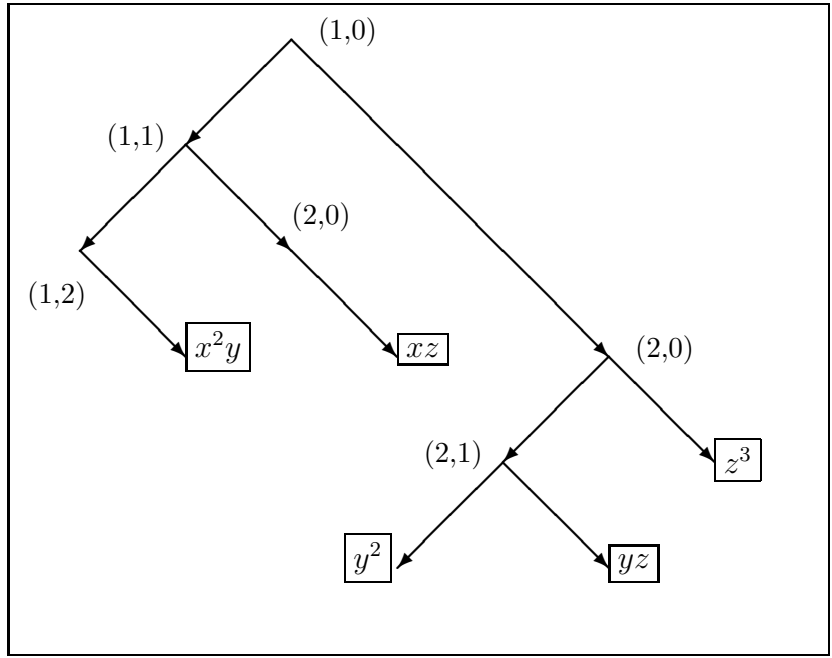
\begin{center}
\fbox{
\begin{picture}(300,240)(20,-40)
\thicklines \put(120,190){\vector(-1,-1){40}}
\put(80,150){\vector(-1,-1){40}}\put(40,110){\vector(1,-1){40}}
\put(80,70){\fbox{$x^2y$}} \put(80,150){\vector(1,-1){40}}
\put(120,110){\vector(1,-1){40}} \put(160,70){\fbox{$xz$}}

\put(120,190){\vector(1,-1){120}} \put(240,70){\vector(-1,-1){40}}
\put(200,30){\vector(1,-1){40}} \put(240,-10){\fbox{$yz$}}
\put(200,30){\vector(-1,-1){40}} \put(140,-10){\fbox{$y^2$}}

\put(240,70){\vector(1,-1){40}} \put(280,26){\fbox{$z^3$}}

\put(130,190){\small (1,0)} \put(50,150){\small (1,1)}
\put(20,90){\small (1,2)} \put(120,120){\small (2,0)}
\put(250,70){\small (2,0)} \put(170,30){\small (2,1)}
\end{picture}}
\end{center}

\noindent
Consider now the structure of Janet tree of the general form as a set
$JT:=\cup \{\nu\}$
of internal nodes and leaves which corresponds to a nonempty binomial set.
To every element $\nu$ of the tree we shall assign the set of five elements
$\nu=\{v,\, d,\,
nd,\, nv,\, nb\}$ with the following structure:
$$
\begin{array}{lcll}
\var(\nu )&=&v& \, \mbox{is the index of the current variable} \\
\dg(\nu )&=&d& \, \mbox{is the degree of the current variable} \\
\ndg(\nu )&=&nd& \, \mbox{is the pointer to the next node in degree} \\
\nvr(\nu )&=&nv& \, \mbox{is the pointer to the next node in variable} \\
\bnm(\nu )&=&bn& \, \mbox{is the pointer to binomial}
\end{array}
$$
In the absence of a child we shall assign the value {\em nil} to
the corresponding pointer. Wherever it does not lead to
misunderstanding we shall identify the pointers $nd$ and $nv$
with the nodes they point out. To the root of $JT$ we assign $\nu_0$
with $\var(\nu_0)=1$ in accordance with labeling (\ref{order}) and
$\dg(\nu_0)=0$.

\noindent
The internal nodes and leaves of tree $JT$ are characterized by the states:
$$
\begin{array}{rl}
\mbox{\bf Internal node:}& ((nv \neq nil \wedge v < \var(nv)) \vee
                 (nd \neq nil \wedge d < \dg(nd)) \\
               &  \phantom{} \wedge bn=nil \\
\mbox{\bf Leaf:}& nv=nil \wedge nd =nil \wedge bn\neq nil
              \wedge d=\dg(\lm(bn)).
\end{array}
$$

For a fast search for Janet divisor in the given tree one can use
the following algorithm {\bf J-divisor} which is an adaptation to the above structure
of Janet tree of the algorithm described in~\cite{GBY1}.

\begin{algorithm}{J-divisor($JT,\,w$)\label{J-divisor}}
\begin{algorithmic}[1]
\INPUT $JT$, a Janet tree; $w$, monomial
\OUTPUT $bn$, a binomial such that $\lm(bn)\, \mid_J \,w$, \\ or $nil$,
 otherwise
 \STATE $\nu:=\nu_0$
 \WHILE{$\deg_{\var(\nu)}(w)\geq \dg(\nu)$}
     \WHILE{$\ndg(\nu)$ {\bf and} $\deg_{\var(\ndg(\nu))}(w)\geq \dg(\ndg(\nu))$}
        \STATE $\nu:=\ndg(\nu)$
     \ENDWHILE
     \IF{$\nvr(\nu)$}
       \STATE $\nu:=\nvr(\nu)$
     \ELSIF{$\ndg(\nu)$}
       \RETURN $nil$
     \ELSE
       \RETURN $\bnm(\nu)$
     \ENDIF
 \ENDWHILE
 \RETURN $nil$
\end{algorithmic}
\end{algorithm}

\noindent
Apparently, the next theorem formulated and proved in~\cite{GBY1} is valid for the
adapted algorithm as well.

\begin{theorem}
 Let $d$ be the maximal total degree of the leading monomials of binomials in $n$
 variables which constitute the finite set $U$. Then the complexity bound of the algorithm
 ${\bf J-divisor}$ and the binary search algorithm is given by
\begin{eqnarray*}
 t_{\bf J-divisor} &=& O(d+n), \\
 t_{\bf BinarySearch} &=& O(n((d+n)\log(d+n)-n\log(n)-d\log(d))).
\end{eqnarray*}
\end{theorem}

\noindent
Thus, the complexity bound for the search of Janet divisor
is $O(n+d)$ where $n$ is the number of variables and $d$ is the maximal
degree of the leading monomials in the binomial basis. Since this bound is even
lower than that for the binary search algorithm, one can expect that the involutive
completion of binomial ideals may be faster than the reduced \Gr basis
completion.

\subsection{Algorithms for Binomial Janet Bases}
Given the generating binomial set $F$ of a toric ideal $\mathcal{I_A}$,
the following algorithm {\bf BinomialJanetBasis} which is a special form
of the general polynomial algorithm~\cite{GB,GBY2} constructs a Janet
basis of $\mathcal{I_A}$.

\begin{algorithm}{BinomialJanetBasis($F,\,\prec$)\label{JanetBasis}}
\begin{algorithmic}[1]
\INPUT $F\in \R\setminus \{0\}$, a finite binomial set; $\prec$, an admissible \\
 \hspace*{0.7cm} ordering
\OUTPUT $G$, a Janet basis of the ideal generated by $F$
\STATE {\bf choose} $f\in F$ with the lowest $\lm(f)$ w.r.t. $\succ$
\STATE $T:=\{f,\lm(f),\emptyset\}$
\STATE $Q:=\{\ \{q,\lm(q),\emptyset\} \mid q\in F\setminus \{f\}\ \}$
\STATE $Q:=${\bf JanetReduce}$(Q,T)$
  \WHILE{$Q \neq \emptyset$}
    \STATE {\bf choose} $p \in Q$ with the lowest $\lm(\bin(p))$ w.r.t. $\succ$
    \STATE $Q:=Q\setminus \{p\}$
    \IF{$\lm(\bin(p)) = \anc(p)$}
          \FORALL{$\{\ r \in T \mid lm(\bin(r))\succ \lm(\bin(p))\ \}$}
          \STATE $Q:=Q \cup \{r\}$; \hspace*{0.4cm} $T:=T \setminus \{r\}$
          \ENDFOR
          \STATE $p:=NF_J(\bin(p),T)$
    \ENDIF
    \STATE $T:=T \cup \{p\}$
    \FORALL{$q\in T$ {\bf and} $x\in NM_J(\bin(q),T)\setminus \nmp(q)$}
      \STATE $Q:=Q \cup \{\ \{\bin(q)\cdot x,\anc(q),\emptyset\}\ \}$
      \STATE $\nmp(q):=\nmp(q)\cap NM_J(\bin(q),T)\cup \{x\}$
    \ENDFOR
    \STATE $Q:=${\bf JanetReduce}$(Q,T)$
  \ENDWHILE
  \RETURN $G:=\{\ \bin(f)\mid f\in T\ \}$
\end{algorithmic}
\end{algorithm}

\noindent
As well as in~\cite{GBY2} to apply the involutive criteria and avoid repeated
prolongations we shall endow with every binomial $f\in F$ the
triple structure
$$
p=\{f,\, u,\, vars\}
$$
such that
$$
\begin{array}{lcl}
\bin(p)&=&f\ \mbox{is binomial itself},\\
\anc(b)&=&u\ \mbox{is the leading monomial of a binomial ancestor of}\ f\
 \mbox{in}\ F\\
\nmp(p)&=&vars \ \mbox{is a (possible empty) subset of variables}.
\end{array}
$$
Here the {\em ancestor} of $f$ is a polynomial $g\in F$ with $u=\lm(g)$
and such that $u\mid \lm(p)$. Moreover, if $\deg(u)<\deg(\lm(p))$, then
every variable occurring in the monomial $\lm(p)/u$ is nonmultiplicative for
$g$. Besides, for the ancestor $g$ the equality $\anc(g)=\lm(g)$ must hold.
These conditions mean that polynomial $p$ was obtained in the course of
the below algorithm {\bf BinomialJanetBasis} from $g$ by a sequence of
nonmultiplicative prolongations.
This tracking of the history in the algorithm allows one to use
the involutive analogues of Buchberger's criteria to avoid unnecessary
reductions.

The set $vars$ contains those nonmultiplicative variables
that have been already used in the algorithm for construction of
nonmultiplicative prolongations. This set serves to prevent repeated
prolongations.

In order to provide minimality of the output Janet basis we
separate~\cite{GB,GBY2} the whole polynomial data into two subsets which
are contained in sets $T$ and $Q$. Set $T$ is a part of
the intermediate binomial basis. Another part of the intermediate basis is contained
in set $Q$ together with all the nonmultiplicative prolongations
of polynomials in $T$ which must be examined in accordance to the
above definition~(\ref{J_basis}) of Janet bases. In so doing, after every
insertion of a new element $p$ in $T$ all elements $r\in T$ such that
$$\lm(\bin(r))\succ \lm(\bin(p))$$ are moved from $T$ to $Q$ as the {\bf for}-loop
6-11 in algorithm {\bf BinomialJanetBasis} does. Such a displacement guaranties
that the output basis is minimal~\cite{GB}.

It should also be noted that for any triple $p\in T$ the set $vars$
must always be a subset of the set of nomultiplicative variables for $\bin(p)$
\begin{equation}
vars\subseteq NM_J(\bin(p),T)\,. \label{cond}
\end{equation}
In the description of algorithm {\bf JanetBinomialBases}
we use the contractions:
\begin{eqnarray*}
&& NM_J(\bin(p),T)\equiv NM_J(\bin(p),\{\bin(f)\mid f\in T\})\,, \\
&& NF_J(\bin(p),T)\equiv NF_J(\bin(p),\{\bin(f)\mid f\in T\})\,, \\
\end{eqnarray*}
\noindent
The insertion of a new polynomial in $T$ may generate new
nonmultiplicative prolongations of elements in $T$ which are added to $Q$ in line 16.
To avoid repeated prolongations the set $\nmp(q)$ of Janet nonmultiplicative
variables for $g$ has been used to construct its prolongations is enlarged with
$x$ in line 17. The intersection placed in this line preserves the condition
(\ref{cond}).

The subalgorithms {\bf JanetReduce} and $\mathbf{NF_J}$ perform Janet reduction of
polynomials in $Q$ modulo polynomials in $T$ and presented below.
In addition to reductions in lines 4 and 19, the Janet normal form computation
is placed in line 12. This is because the replacement of elements from $T$ to $Q$ may lead
to the tail reducibility of the binomial in $p$. Such a reducibility may be caused
by converting of some nonmultiplicative variables for binomials in $T$ into multiplicative
due to the replacement.

In subalgorithm {\bf JanetReduce} computation of the Janet normal
form $h$ is done in line 6 for every binomial $\bin(p)$ in $T$. If
$h$ is nonzero, then line 8 checks if $\lm(\bin(p))$ was subjected
by reduction. If the reduction took place $\lm(h)$ cannot be
multiple of any monomial in the set $\{\lm(\bin(g)) \mid g\in
T\}$~\cite{GB}. Therefore, one has to insert the triple with $h$
in the output set $Q$ as shown in line 9 as $h$ cannot have
ancestors among polynomials in $T$ and one must also examine all
nonmultiplicative prolongations of $h$. If $\lm(\bin(p))$ is Janet
irreducible modulo $\{\lm(\bin(g)) \mid g\in T\ \}$, then the
triple $\{h,\anc(p),\nmp(p)\}$ is added to $Q$ in line 11.

\begin{algorithm}{JanetReduce($Q,T$)}
\begin{algorithmic}[1]
\INPUT $Q$ and $T$, sets of triples
\OUTPUT $Q$ whose polynomials are Janet head reduced modulo $T$
\STATE $S:=Q$
\STATE $Q:=\emptyset$
\WHILE {$S \neq \emptyset$}
   \STATE {\bf choose} $p\in S$
   \STATE $S:=S\setminus \{p\}$
   \STATE $h:=\mathbf{NF_J}(p,T)$
   \IF{$h\neq 0$}
       \IF{$\lm(\bin(p))\neq \lm(h)$}
       \STATE $Q:=Q\cup \{h,\lm(h),\emptyset\}$
       \ELSE
       \STATE $Q:=Q\cup \{h,\anc(p),\nmp(p)\}$
       \ENDIF
   \ENDIF
\ENDWHILE
\RETURN $Q$
\end{algorithmic}
\end{algorithm}

\noindent
Subalgorithm $\mathbf{NF_J}(p,T)$ performs the Janet reduction of a
binomial $g=\bin(p)$ modulo polynomial set in $T$:

\begin{algorithm}{$\mathbf{NF_J}(f,T)$}
\begin{algorithmic}[1]
\INPUT $f=\{\bin(f),\anc(f),\nmp(f)\}$, a triple; $T$, a set of triples
\OUTPUT $h=NF_J(\bin(f),T)$, the Janet normal form of the \\
  \hspace*{1.2cm}binomial in $f$ modulo binomial set in $T$
    \STATE $G:=\{\bin(g)\mid g\in T\}$
    \STATE $h:=\bin(f)$
    \IF{$\lm(h)$ is Janet reducible modulo $G$}
       \STATE {\bf choose} $g\in T$ such that $\lm(\bin(g))\mid_J \lm(h)$
       \IF{$\lm(h) \neq \anc(f)$ {\bf and} \\
       \hspace*{0.4cm}{\bf CriterionI}$(f,g)$ {\bf or} {\bf CriterionII}$(f,g)$}
            \RETURN $0$
       \ENDIF
    \ELSE
       \WHILE{$h\neq 0$ {\bf and} $h$ has a term $t$ Janet reducible modulo $G$}
            \STATE {\bf choose} $q\in G$ such that $\lm(q)\mid_J t$
            \STATE $h:=h - q\cdot t/\lm(q)$
       \ENDWHILE
    \ENDIF
   \RETURN $h$
\end{algorithmic}
\end{algorithm}
For the head reducible input binomial $\bin(f)$ the two criteria are
verified in line 5:
\begin{itemize}
\item {\bf CriterionI}$(f,g)$ is true iff $\anc(f)\cdot \anc(g) \mid \lm(\bin(f))$.
\item {\bf CriterionII}$(f,g)$ is true iff
$\deg(\lcm(\anc(f)\cdot \anc(g)))< \deg(\lm(\bin(f))$.
\end{itemize}
These criteria are the Buchberger criteria~\cite{Buch85}
adapted to the involutive completion procedure. If any of the two criteria is
true, then $NF(\bin(f),T)=0$~\cite{GBY2}.

It should be noted that the Janet normal form
is uniquely defined and, hence, uniquely computed by the above subalgorithm.
This uniqueness hold because of the uniqueness of a Janet divisor among the leading terms of binomials
in $T$ at every step of intermediate computations~\cite{GB}.

\section{Examples}

As we emphasized in Sect.3.1, in the course of involutive completion of the initial binomial
generators for a toric ideal the reduction can be performed very fast due to the fast search
for a Janet divisor, This fast search is provided by the use of the Janet tree structures
for intermediate binomial set. Our computer experiments with C/C++ codes implementing
polynomial algorithms for Janet bases~\cite{GBY2} perfectly strengthen this theoretical
fact. In particular this fast reduction in addition to suppressing swell of intermediate
integer coefficients results in high computational speed observed for the benchmark
collection used for testing \Gr bases software~\cite{GBY2}. These benchmarks, however,
are not very "sparse" with respect to degrees of variables occurring in the generating set.
By contrast, the generating binomial sets for toric ideals especially for those arising in integer
programming problem are usually highly sparse. This may lead to much larger cardinality
of a Janet basis than that of the reduced \Gr basis and thereby annihilate the advantages of
involutive reduction.

Consider the example taken from~\cite{Morales}
$$ \mathcal{I_A}=\{\ x_0x_1x_2x_3x_4-1, x_2^{29}x_3^5-x_1^{14}x_4^{20}, x_1^{39}-x_2^{25}x_3^{14}\ \}\,.$$
Our C++ package~\cite{GBY2} generates the degree-reverse-lexicographical Janet basis
of $\mathcal{I_A}$ with 7769 binomials whose sorting with respect to the ordering chosen
gives
$$ \{\ x_0x_1^3x_3x_4^{281} - x_1x_2^{280}, x_0x_2^{61}x_3^2x_4^{221} - x_1x_2^{279},
x_0x_1^2x_3x_4^{281} - x_2^{280},\ldots,x_0x_1x_2x_3x_4-1\ \}\,
$$
where we explicitly show only three highest ranking binomials and the lowest one. The computing time
on a Pentium III 700 Mhz based PC running under RedHat Linux 6.2 is 6 seconds that is noticeably larger
than the running time for direct computation of the reduced \Gr basis which contains 19 binomials only:
\begin{eqnarray*}
&&\{\ x_{{0}}{x_{{1}}}^{2}x_{{3}}{x_{{4}}}^{281}-{x_{{2}}}^{280},
   {x_{{2}}}^{281}-x_{{1}}{x_{{4}}}^{280},x_{{0}}{x_{{3}}}^{2}{x_{{4}}}^{221}-x_{{1}}{x_{{2}}}^{218},
   {x_{{1}}}^{2}{x_{{2}}}^{219}-x_{{3}}{x_{{4}}}^{220}, \\
&&\ x_{{0}}{x_{{3}}}^{3}{x_{{4}}}^{161}-{x_{{1}}}^{4}{x_{{2}}}^{156},
   {x_{{1}}}^{5}{x_{{2}}}^{157}-{x_{{3}}}^{2}{x_{{4}}}^{160},
   x_{{0}}{x_{{3}}}^{4}{x_{{4}}}^{101}-{x_{{1}}}^{7}{x_{{2}}}^{94},
   {x_{{1}}}^{8}{x_{{2}}}^{95}-{x_{{3}}}^{3}{x_{{4}}}^{100},\\
&&\ x_{{0}}{x_{{1}}}^{4}{x_{{4}}}^{61}-{x_{{2}}}^{61},
   {x_{{2}}}^{62}x_{{3}}-{x_{{1}}}^{3}{x_{{4}}}^{60},
   x_{{0}}{x_{{3}}}^{5}{x_{{4}}}^{41}-{x_{{1}}}^{10}{x_{{2}}}^{32},
   {x_{{1}}}^{11}{x_{{2}}}^{33}-{x_{{3}}}^{4}{x_{{4}}}^{40},\\
&&\ x_{{0}}{x_{{2}}}^{26}{x_{{3}}}^{15}x_{{4}}-{x_{{1}}}^{38},
   {x_{{1}}}^{39}-{x_{{2}}}^{25}{x_{{3}}}^{14},
   x_{{0}}{x_{{1}}}^{15}{x_{{4}}}^{21}-{x_{{2}}}^{28}{x_{{3}}}^{4},
   {x_{{2}}}^{29}{x_{{3}}}^{5}-{x_{{1}}}^{14}{x_{{4}}}^{20},\\
&&\ x_{{0}}{x_{{3}} }^{10}{x_{{4}}}^{21}-{x_{{1}}}^{24}{x_{{2}}}^{3},
    {x_{{1}}}^{25}{x_{{2}}}^{4}-{x_{{3}}}^{9}{x_{{4}}}^{20},x_{{0}}x_{{1}}x_{{2}}x_{{3}}x_{{4}}-1\ \}\,.
\end{eqnarray*}
Accordingly, such a computer algebra system as {\em Singular}~\cite{Singular} needs
much less than 1 second to compute this \Gr basis on the same computer.

\noindent
Having ascertained this drawback of the involutive method with respect to the \Gr basis one in computing
toric ideals we designed another algorithmic approach to computing \Gr bases~\cite{GB02}.
This approach preserves the Janet-like tree structure and uniqueness of a divisor though underlying
division is not involutive since it does not satisfy the axioms in~\cite{GB}. On the other hand
the resulting bases unlike Janet bases are often reduced as \Gr bases and their cardinality
is always less or equal to the cardinality of Janet bases. For toric ideals the new bases are much more
compact then Janet bases. We have not implemented yet the new algorithm and so we demonstrate
the compactness of its output in comparison with algorithm {\bf BinomialJanetBasis} by the following
simple example taken from~\cite{BLR99}:
$$\mathcal{I_A}=\{\ x^7-y^2z, x^4w-y^3, x^3y-zw\ \}\,.$$
The reduced \Gr basis and Janet basis of this toric ideal for the degree-reverse-lexico\-graphic
order induced by $x\succ y\succ z\succ w$ are
$$ \{\ x^7-y^2z, x^4w-y^3, x^3y-zw, y^4-xzw^2\ \}$$
and
\begin{eqnarray*}
&& \{\ x^7-y^2z, x^6y-x^3zw, x^6w-x^2y^3, x^5y-x^2zw, x^2y^4-x^3zw^2, x^5w-xy^3, \\
&&\ \  x^4y-xzw, x^2zw^2-xy^4, x^4w-y^3, x^3y-zw, y^4-xzw^2\ \}\,,
\end{eqnarray*}
respectively. Their cardinalities are 4 and 11. The new basis contains 5 elements
$$ \{\ x^7-y^2z, x^4y-xzw, x^4w-y^3, x^3y-zw, y^4-xzw^2\ \}$$
and contains only single extra element in comparison  with the reduced \Gr basis.

It should be noted that there are also a number of other efficient algorithms
computing \Gr bases of toric ideals (see, for example,~\cite{BLR99,DBUr95,HoSt95})
which are differ greatly from just completion of a generating binomial set to a \Gr basis.
After implementation of our new algorithm we are planning to run the underlying code
for collection of large examples given in~\cite{Sturmfels96,BLR99} and other references.

\section{Acknowledgements}

The work was supported in part by the RFBR grants 00-15-96691, 01-01-00708 and
by grant Intas 99-1222.


\begin{thebibliography}{99}

\bibitem{Sturmfels96} Sturmfels, B.: {\em \Gr bases and convex
 polytopes.} University Lecture Series {\bf 8}, Providence, RI, American Mathematical Society, 1996.

\bibitem{BLR99} Bigatti, A.M., La Scala, R., Robbiano, L.:
Computing toric ideals.  {\em J. Symb. Comp.} {\bf 27} (1999) 351-365.

\bibitem{KM96} Koppenhagen, U., Mayr, E.W.: An Optimal Algorithm for
Constructing the Reduced \Gr Basis of Binomial Ideals. {\em Proceedings of ISSAC'96},
ACM Press, 1996, 55-62.

\bibitem{AL94} Adams, W.W., Loustaunau, P.: {\em An Introduction to \Gr Bases}.
Graduate Studies in Mathematics {\bf 3}, American Mathematical Society,
1994.

\bibitem{Pot94} Pottier, L.:
\Gr Bases of Toric Ideals. {\em Rapport de recherche} {\bf 2224} (1997),
INRIA Sophia Antipolis.

\bibitem{GB} Gerdt, V.P., Blinkov, Yu.A.: Involutive Bases of Polynomial
 Ideals. {\em Math. Comp. Simul.} {\bf 45} (1998) 519-542; Minimal Involutive
 Bases. {\em Math. Comp. Simul.} {\bf 45} (1998) 543-560.

\bibitem{GBY1} Gerdt V.P., Blinkov Yu.A., Yanovich D.A.
 Construction of Janet Bases. I. Monomial Bases. In: {\em Computer Algebra in
 Scientific Computing / CASC'01},
 V.G.Ganzha, E.W.Mayr and E.V.Vorozhtsov (Eds.), Springer-Verlag,
 Berlin, 2001, pp.233-247.

\bibitem{GBY2} Gerdt, V.P., Blinkov, Yu.A., Yanovich, D.A.:
 Construction of Janet Bases. II. Polynomial Bases. In: {\em Computer Algebra in
 Scientific Computing / CASC'01},
 V.G.Ganzha, E.W.Mayr and E.V.Vorozhtsov (Eds.), Springer-Verlag,
 Berlin, 2001, pp.249-263.

\bibitem{CoTr91} Conti, P., Traverso, C.:
Buchberger algorithm and integer programming. {\em Proceedings of AAECC-9},
Springer LNCS {\bf 539} (1991) 130-139.

\bibitem{B01} Blinkov, Yu.A.: Method of Separative Monomials
for Involutive Divisions. {\em Programming and Computer Software}
{\bf 3} (2001) 43-45.

\bibitem{Janet} Janet, M.: {\em Le\c cons sur les Syst\`emes
 d'Equations aux D\'eriv\'ees Partielles}, Cahiers Scientifiques {\bf IV},
 Gauthier-Villars, Paris, 1929.

\bibitem{Buch85} Buchberger, B.: Gr\"obner Bases: an Algorithmic
 Method in Polynomial Ideal Theory. In: {\em Recent Trends in
 Multidimensional System Theory}, N.K. Bose (ed.), Reidel, Dordrecht
 (1985) 184--232.

\bibitem{Morales} Morales, M.: Equations des vari\'et\'es monomiales. {\em Preprint},
 Universit\'e de Grenoble I, France.

\bibitem{GB02} Gerdt, V.P., Blinkov, A.Yu.: Involutive-like \Gr Bases. {\em In preparation}.

\bibitem{DBUr95} Di Biase, F., Urbanke, R.:
An algorithm to calculate the kernel of certain polynomial ring
homomorphisms. {\em Experimental Mathematics} {\bf 4} (1995) 227-234.

\bibitem{HoSt95} Hosten, S., Sturmfels, B.:
GRIN: An implementation of Groebner basis for integer programming.
In: {\em Integer Programming and Combinatorial Optimization}, Balas, E., Clausen J.,
(Eds.), LNCS {\bf 920}, Springer-Verlag, New York, 1995, pp.267-276.

\bibitem{Singular} Greuel, G.-M., Pfister, G., Sch\"{o}enemann, H.:
{\em Singular: A Computer Algebra
 System for Polynomial Computation}, Department of Mathematics,
 University of Keiserslautern (2001)
 http://www.singular.uni-kl.de/Manual/2-0-0/.


\end{thebibliography}
\end{document}